\documentclass[12pt]{amsart}
\usepackage{amsmath,amssymb}
\usepackage{amsfonts}
\usepackage{amsthm}
\usepackage{latexsym}
\usepackage{graphicx}
\usepackage{epsfig}
\usepackage{setspace}
\usepackage{xcolor}
\newtheorem{prop}{Proposition}[section]
\newtheorem{remark}[prop]{Remark}

\makeatletter \@addtoreset{equation}{section} \makeatother

\def\ppt{\frac{\partial}{\partial t}}

\def\RR{{\mathrm R}}
\def\WW{{\mathrm W}}

\def\Rc{{\mathrm {Rc}}}

\def\SS{{\mathrm S}}

\def\He{\mathrm {Hess}}
\def\lie{\mathcal{L}}

\begin{document}

\title[Dedicated to Richard S. Hamilton]{Geometry and Analysis of Gradient Ricci Solitons in Dimension Four}
\dedicatory{Dedicated to Richard S. Hamilton on forty years of Ricci flow.}
\author{Xiaodong Cao$^*$
}

\address{Department of Mathematics,
 Cornell University, Ithaca, NY 14853-4201}
\email{xiaodongcao@cornell.edu}
%\curraddr{}
%\thanks{}

\author{Hung Tran$^{\flat}$}
%\thanks{$^{\flat}$Research
%partially supported by NSF grant no. \# }

\address{Department of Mathematics,
 Texas Tech University}
\email{hung.tran@ttu.edu}
%\curraddr{}
%\thanks{}

%    General info
%\renewcommand{\subjclassname}{%
 % \textup{2000} Mathematics Subject Classification}
%\subjclass[2000]{Primary 53C44}
% Global differential geometry
% 53C21 Methods of Riemannian geometry, including PDE methods;
%   curvature restrictions [See also 58J60]
% 53C44 Geometric evolution equations (mean curvature flow)
% 53C55 Hermitian and K\"ahlerian manifolds [See also 32Cxx]
% Qualitative properties of solutions
% 35B35 Stability, boundedness
%Parabolic equations and systems [See also 35Bxx, 35Dxx, 35R30,
%35R35, 58J35]
%35K55 Nonlinear PDE of parabolic type
%35K90 Abstract parabolic evolution equations
% Partial differential equations on manifolds; differential operators
%58J35 Heat and other parabolic equation methods
% 58J37 Perturbations; asymptotics
% Systems theory; control - Stability
% 93D05 Lyapunov and other classical stabilities
%       (Lagrange, Poisson, $L^p, l^p$, etc.)

\date{\today}

\maketitle

%{\centering\footnotesize Dedicated to Richard S. Hamilton on forty years of Ricci flow.\par}

\begin{abstract} Gradient Ricci solitons have garnered significant attention both as self-similar solutions and singularity models of the Ricci flow. This survey article  starts with a list of %{%non-trivial} 
examples; it also provides some geometric aspects of gradient Ricci solitons, including various asymptotic behaviors;   finally, it discusses some recent results on classification and rigidity. In particular, this survey focuses on dimension four.
\end{abstract}

%%%%%%%%%%%%%%%%%%%%%%%%%%%%%%%%%%%%
\section{Introduction}
%09/19/24-Hung

In his seminal paper \cite{H3}, R. Hamilton introduced the Ricci flow equation
\begin{equation}
 \label{rf}
 \ppt g= -2 \Rc, \end{equation} aiming to find the ``Best" metric on a manifold $M$.
%Talk about RF here, Poincare, Bohm-Wilking, Brendle-Schoen, etc.
Hamilton established the foundation and fundamental theory via a series of pioneering articles including, but not limited to, \cite{H3, HPCO, H93harnack, H93eternal, Hsurvey}. His deep understanding and breath of knowledge led to a vision that inspired several generations of mathematicians. 

In particular, the program has been carried out by G. Perelman in his proof of the Poincar\'e and Thurston's Geometrization conjectures (cf. \cite{perelman1, perelman2, perelman3}). Later, the Ricci flow method was developed further  by C. B\"ohm and B. Wilking in the proofs of the Space-form Theorem for  manifolds with 2-positive curvature operators  \cite{bohmwilking}; and by S. Brendle and R. Schoen in their proof of the Differentiable-Sphere Theorem \cite{  bs091, bs072}. The Ricci flow machinery, combining topology, geometry, and analysis in a single heat-type equation, proves to be a powerful and flexible tool. 

Currently, the study of Ricci flows is a tremendously vast and active area of research. There are several excellent surveys on the subjects such as \cite{CCZ08, caohd09,bamlersurvey21, CKM24}. In the present article, we restrict ourselves to gradient Ricci solitons, aka singularity models of the Ricci flow, and give an overview of our perspective on its state of the art, with a particular emphasis on dimension four. We acknowledge that there are many other important concerns not addressed in this survey, hence by no means it is comprehensive.

%Set up notation here.
\subsection{Notation and Generalities}
Throughout this paper, $(M^n, g)$ is a Riemannian manifold of dimension $n$. $\RR, \WW, \Rc$, $\text{K}$ and $\SS$ stand for the Riemannian curvature operator, Weyl tensor, Ricci curvature, sectional curvature and scalar curvature respectively. For a finite dimensional real vector space (bundle) $V$, $\Lambda^{2}(V)$ denotes the space of bivectors or two-forms. In our case,  the space of interest is normally the tangent bundle and when the context is clear, the dependence on $V$ is omitted. 
    
Given an orthonormal basis  $\{E_{i}\}_{i=1}^{n}$ of $T_p M$, it is well-known that we can construct an  orthonormal frame around p such that $e_i(p)=E_i$ and $\nabla e_i \mid_p =0$. Such a frame is called normal at p; let $e^{i}$ be the dual one-form via the Riemannian metric $g$ and $e_{ij}$ be the shorthand notation for $e^{i}\wedge e^{j}\in \Lambda^{2}$. The Einstein summation convention is also assumed throughout. Last but not least, $\delta$, $\Delta$, and $\Delta_f$ denote the divergence, Laplacian, and drift Laplacian. Our convention agrees with \cite{pe06book}. 

 %Peterson'text Exercise 5 page 56 

%The drifted Laplacian is defined as 
%\begin{equation*} \Delta_{f} =\Delta-\nabla_{\nabla f},
%\end{equation*}
%where $f$ is a smooth function on $M$. With respect to a normal frame, for any $(m,0)$-tensor $T$, the divergence operator is defined as
%\[(\delta T)_{p_{2}...p_{m}}=\sum_{i} \nabla_{i}T_{i p_{2}...p_{m}}=\nabla_{i}T^i_{p_{2}...p_{m}}; \]
%while its interior product by a vector field $X$ is defined as
%\[(i_{X} T)_{p_{2}...p_{m}}=T_{X p_{2}...p_{m}}. \]

Equation (\ref{rf}) is only weakly parabolic, it will generally produce finite-time singularities, which lead to the study of singularity models. It is one of the most central problems in this field to understand and classify these models. Probably the most important singularity model is the Ricci soliton, which is a self-similar solution to the Ricci flow equation (\ref{rf}), and arises as a finite-time singularity model. Algebraically, it can be characterized by the equation
\begin{equation}
 \label{rs}
 \Rc+\frac{1}{2}\lie_{X}~g=\lambda g, \end{equation}
where $\lie$ is the Lie derivative and $\lambda$ is a real number. In the special case that the vector field $X$ is a gradient vector field $\nabla f$, this is called gradient Ricci soliton (GRS),
\begin{equation}\label{grs}
	\text{Rc}+\frac{1}{2}\lie_{\nabla f} ~g =\Rc+ \He~ f= \lambda g.
\end{equation}
Hence, $f$ is called the potential function; while $\lambda =+1,~ 0,~ -1$ corresponding to the case of shrinking, steady and expanding solitons.

 The following formulas are direct consequence of (\ref{grs}) by simple algebraic manipulation and the Bianchi identities (for example, see \cite{chowluni}):
\begin{align}
\SS + \triangle f &= n\lambda,\\
\label{rcdiv}
\frac{1}{2}\nabla_{i}\SS=\nabla^{j}\RR_{ij} &= \RR_{ij}\nabla^{j}f.
%\label{rcandf}
\end{align}
Due to its symmetry, the Ricci curvature is frequently considered as an endomorphism on $TM$. Thus, via the second Bianchi's identity, we can deduce that
\begin{equation}
	\label{rcandf}
	\Rc(\nabla{f})=\frac{1}{2}\nabla{\SS}=\delta \Rc.
\end{equation}
%Recall that, for a normal frame,
%\[ \delta \Rc (X)= \sum_{i} g((\nabla_{e_i}\Rc)X, e_i).\]
Consequently, the following is considered as a conservation law,
\begin{equation}
	\label{nablafandS}
	\SS+|\nabla f|^2-2\lambda f = \text{constant}.	
\end{equation}
Here is another interesting identity \cite{chowluni}, 
\begin{align}
	\label{lapS}
	\triangle\SS+2|\text{Rc}|^2 &= \left\langle{\nabla f,\nabla \SS}\right\rangle+2\lambda \SS.
\end{align}

What is more, equation (\ref{grs}) can be considered as a generalization of the Einstein equation. The study of the Bakry-Emery Ricci tensor $ \Rc+\He{f}$ draws independent interest from the areas of probability and classical Riemannian geometry, which we will not address in this survey.
Similar to the Einstein structure, the Ricci soliton structure automatically carries some restriction onto its topology and geometry.  If $\lambda \geq 0$, then the maximum principle and equation $(\ref{lapS})$ imply that 
\[\SS\geq 0.\] 
Moreover, such a complete GRS has positive scalar curvature unless it is isometric to the flat Euclidean space  \cite{zhang09completeness, chenbl09}. In case $\lambda<0$, Z. Zhang \cite{zhang09completeness} shows that $\SS$ is also bounded below.

On the topological side, W. Wylie \cite{wylie07} shows that any complete shrinking GRS has finite fundamental group, generalizing previous work of \cite{derdzinski06, fg08}. Indeed, the proof is applicable for a larger class of Riemannian manifolds with a positive lower bound on the Bakry-Emery Ricci tensor. For a compact manifold, a finite fundamental group leads to the vanishing of the first Betti number \cite{derdzinski06}. %(and the $(n-1)$-th, due to Poincare's duality). 
Later, the result is improved by B. Chow and P. Lu \cite{cl16}, determining what the finiteness depends on. 

In \cite{MW17compact}, O. Munteanu and J. Wang prove that a shrinking Ricci soliton with nonnegative sectional curvature and positive Ricci curvature must be compact. Their proof is by contradiction based on a clever estimate on a suitable function constructed using the smallest eigenvalue of the Ricci curvature and an integral inequality derived by H.-D. Cao and D. Zhou \cite{caozhou10}. %[see J. Differential Geom. 85 (2010), no. 2, 175–185; MR2732975]

% that the scalar curvature S must satisfy, outside a geodesic ball, $S>b\log f$ for some constant $0<b\leq 1.$ Combining this estimate with the integral estimate of , the authors conclude that M must be compact.

%XC: Should we mention about their proof here? I think so to make the paragraph more meaningful

%	Fang-Man-Zhang 2008: complete GRS has finite topological type. (XC: I am not sure if we should include this or not, I prefer not)
	In the complex setup, it is natural to consider a K\"{a}hler gradient Ricci soliton (KGRS) $(M, g, f, J, \lambda)$. That is, for an almost complex structure $J$, $(M, J, g)$ is K\"{a}hler and $(M, g, f, \lambda)$ is a GRS. KGRS arise naturally in the context of running a Ricci flow preserving a K\"{a}hler structure. Thus, this topic has extensive literature; see, for examples, \cite{tian1, WZ04toric, caohd09, MW15topo, CS2018classification, CF16conical, CD2020expanding, DZ2020rigidity} and references therein. We will only mention results in the K\"{a}hler case briefly  in this survey.

The organization of the rest of the article is as follows. In Section 2, we will display some known examples; in Section 3, we will discuss various asymptotic behaviors; in Section 4, we will focus on classification and rigidity results. \\

\textbf{Acknowledgment.} The first author would like to express his deepest gratitude to Professor Richard Hamilton for invaluable mentorship, guidance and unwavering support throughout his career. The authors would like to thank Huai-Dong Cao, Eric Chen, Xiaolong Li, Peng Lu, Ernani Ribeiro Jr. and Detang Zhou for their helpful suggestions. They would also like to thank the anonymous referees for many valuable  comments, including pointing out a few extra references. Cao's research is partially supported by Simons Foundation [585201]; Tran's research is partially supported by grants from the Simons Foundation [709791], the National Science Foundation [DMS-2104988], and the Vietnam Institute for Advanced Study in Mathematics.

%%%%%%%%%%%%%%%%%%%%%%%%%%%%%%%%%%%%%
\section{\textbf{Examples of Ricci Solitons}}%%%%%%%%%%%%%%%
%%%%%%%%%%%%%%%%%%%%%%%%%%%%%%%%%%%%%
In this section, we review examples of GRS. By default, all Einstein manifolds $(N, g_{E})$ with $\Rc(g_E)=\lambda g_E$ are GRS. That is, the potential function $f$ is a constant in space and $\lambda$ is the Einstein constant. Another example is the Gaussian soliton which is essentially the Euclidean space $\mathbb{R}^{n}$ with its standard metric $g_{Euc}$; here the potential function $f=\lambda \frac{|x|^2}{2}$. It is interesting to note that $\lambda$ can assume any arbitrary real value and, thus, the Gaussian soliton can be either shrinking, steady, or expanding. 

A soliton is called {\bf rigid} if it is isometric to a quotient $N\times_\Gamma \mathbb{R}^k$ \cite{PW09rigid}, for an Einstein manifold $N$. Here $\Gamma$ acts freely on $N$ and by orthogonal transformations on $\mathbb{R}^k$ (no translational components), and $f=\frac{\lambda}{2} |x|^2$, where $x$ is the standard coordinate of the Euclidean factor. For $\lambda>0$, such a GRS has a finite fundamental group; thus, it is rigid if and only if it is isometric to a finite quotient $N\times_\Gamma \mathbb{R}^k$. In subsection 4.4, we will discuss a different type of rigidity.

%a product of those above, first observed by P. Petersen and W. Wylie \cite{pewy09}, is called a \textit{rank-k rigid GRS}, namely a quotient of $(N\times \mathbb{R}^k, g_E\times g_{sd})$. %POTENTIAL function?to get a at vector bundle over a base that is Einstein and

Other nontrivial examples of GRS are rare and many of them are K\"ahler. In \cite{tian1}, G. Tian and X. Zhu prove that a compact K\"ahler manifold admits at most one soliton structure up to automorphisms, there are examples on which there is none \cite{WZ04toric}. We summarize a few important constructions in the following table.

\begin{center}
	\begin{tabular}{ |c|c|c|c| } 
		\hline
		Topology & Geometry & Type & Reference \\
		\hline 
		$\mathbb{C}^n$ & Toric symmetry & Non-compact, Steady & \cite{AC23} \\ 
		\hline 
		$\mathbb{C}^n$ & $U(n)$-symmetry and $\text{K}>0$ & Non-compact, Steady & \cite{caohd96} \\ 
		\hline
		$\mathbb{R}^n$ & $SO(n)$-symmetry and $\RR>0$ & Non-compact, Steady & \cite{bryantunpublished, chowluni} \\ 
		\hline
		$M \times \mathbb{R}^{k+1}$ & Doubly warped product & Non-compact, Steady & \cite{Ivey94} \\ 
		\hline
		$X_n$ & $U(n)$-symmetry and $\text{K}\geq 0$ & Non-compact, Steady & \cite{caohd96, CV96} \\
		 \hline 
		 $\mathbb{C}^n$ & $U(n)$-symmetry with cone-like end & Non-compact, Expanding & \cite{caohd97limits, CV96} \\ 
		 \hline
		 $\mathbb{R}^n$ & $SO(n)$-symmetry and $\RR>0$ & Non-compact, Expanding & \cite{chowetc1} \\ 
		 \hline
		$N$ & Homogeneous & Non-compact, Expanding & \cite{Lauret01} \\ 
		 \hline
		 $L^{-k}$,  $k>n$ & $U(n)$-symmetry and  cone-like end & Non-Compact, Expanding & \cite{CV96, FIK03Kahler} \\
		  \hline
		 $M_k$ & $U(n)$-symmetry and $\Rc> 0$ & Compact, Shrinking & \cite{caohd96, koi90} \\
		\hline
		$L^{-k}$,  $0<k<n$ & $U(n)$-symmetry and  cone-like end & Non-Compact, Shrinking & \cite{FIK03Kahler} \\
		\hline
		$\mathbb{CP}^2 \# 2\overline{\mathbb{CP}^2}$ & Toric symmetry & Compact, Shrinking & \cite{tian1, WZ04toric} \\
		\hline
		$\text{Bl}_1(\mathbb{P}^1\times \mathbb{C})$ & Toric symmetry & Non-Compact, Shrinking & \cite{BCCD22KahlerRicci} \\
		\hline
	\end{tabular}
\end{center}

We comment on and further clarify some notations in the table: 
\begin{itemize}
%	\item $K$: sectional curvature, $\RR$: curvature operator
	\item The constructions of steady soliton metrics on $\mathbb{C}^n$ by H.-D. Cao \cite{caohd96}, on $\mathbb{R}^n$ by R. Bryant \cite{bryantunpublished}, and on a doubly warped product by T. Ivey \cite{Ivey94}, all generalize the cigar soliton on $\mathbb{R}^2 \cong \mathbb{C}^1$ discovered by R. Hamilton \cite{Hsurface}. 
	\item $X_n$ is the total space of a line bundle over the complex projective space $\mathbb{CP}^{n-1}$, for $n\geq 2$. In real dimension four, $X_2$ is the tangent bundle to the sphere $\mathbb{CP}^1 \cong \mathbb{S}^2$ (the blow-up of $\mathbb{C}^2/\mathbb{Z}_2$ at the origin). 
	\item $N$ denotes a simply-connected Lie group with a left-invariant metric satisfying certain criteria. Some examples include generalized Heisenberg groups and two-step nilpotent Lie groups \cite{chowetc1}.   
	\item $M_k$ is the total space of the twisted projective line bundle $$\pi: \mathbb{P}(L^k\bigotimes L^{-k})\mapsto \mathbb{CP}^{n-1}$$ for each $n\geq 2$ and $1\leq k\leq n-1$, $L$ is the hyperplane line bundle over $\mathbb{CP}^{n-1}$. In real dimension four, this shows the existence of shrinking  K\"ahler GRS metric on  $\text{Bl}_1\mathbb{CP}^2=\mathbb{CP}^2 \# \overline{\mathbb{CP}^2}$ (the blow-up of $\mathbb{CP}^2$ at one point), here $ \overline{\mathbb{CP}^2}$ denotes taking the opposite orientation of $ \mathbb{CP}^2$.   
	\item $L^{-k}$ is the total space of a twisted line bundle over $\mathbb{CP}^{n-1}$ where $(-k)$ is the twisting number. Thus, $L^{-k}$ is obtained by gluing a copy of $\mathbb{CP}^{n-1}$ into $(\mathbb{C}^n/\{0\})/\mathbb{Z}_k$. In particular, $L^{-1}$ is the blow-up of $\mathbb{C}^n$ at one point. In real dimension four, this shows the existence of shrinking K\"ahler GRS on $\text{Bl}_1\mathbb{C}^2$ (the blow-up of $\mathbb{C}^2$ at one point). 
	\item In real dimension four, on $\text{Bl}_2\mathbb{P}^2\cong \mathbb{P}^2 \# 2\overline{\mathbb{P}^2}$ (the blow-up of $\mathbb{CP}^2$ at two points), the existence of a K\"ahler soliton metric was obtained by X.-J. Wang and X. Zhu \cite{WZ04toric} by solving a Monge-Ampere equation. 
\end{itemize}

There are several extensions of the aforementioned constructions such as \cite{GK04, dw09, DW09expanding, FW11sasaki, dw11, yang12, BDGW15, BDW15, CS2018classification, Jab15, IS17, CS2018classification}. In particular, in a series of papers, A. Dancer, M. Wang, and their coauthors develop a systematic study of solitons structures on multiple warped products. For instance, in \cite{dw11}, A. Dancer and M. Wang develop a framework using group actions with an one-dimensional orbit space; the examples in \cite{caohd96, CV96, FIK03Kahler, koi90} all appear as special cases. 
See also \cite{appleton22, wink23, NW2023} for additional cohomogeneity one expanding and steady solitons of a similar flavor. The work of J. Lauret \cite{Lauret01} is generalized by M. Jablonski \cite{Jab15}. The result of \cite{WZ04toric} is generalized by \cite{CS2018classification, IS17} for Fano 3-folds with torus actions and {G}orenstein del {P}ezzo surfaces. Moreover, it is a common strategy to obtain existence of asymptotically conical expanders under certain conditions (K\"{a}hler by \cite{CD2020expanding}, non-negative curvature operator \cite{SS13, deruelle16cone}, or non-negative scalar curvature \cite{BC2023degree}). For several of these constructions, the metric can be written down explicitly, see \cite{BHJM15, torres18, maximo14}. 
%%%%%%%%%%%%%%%%%%%%%%%%%%%%%%%%%
\section{Asymptotic Behavior}
In this section, we focus on studying the asymptotic behavior of a complete non-compact GRS. We use $r(x)$ to denote the distance to a fixed reference point $x_0$.
\subsection{Estimates on the Potential Function $f$}%%%%%%%%
%%%%%%%%%%%%%%%%%%%%%%%%%%%%%%%%%
The potential function $f$ characterizes the difference between a GRS and an Einstein manifold.  For a shrinking GRS with the normalization $\lambda=\frac{1}{2}$, H.-D. Cao and D. Zhou \cite{caozhou10} prove that
\begin{equation}
	\label{estimatef1}
	\frac{1}{4}(r(x)-c_1)^2\leq f(x) \leq  \frac{1}{4}(r(x)+c_2)^2,
\end{equation}
where $c_1$ and $c_2$ are positive constants depending only on the dimension $n$ and the geometry of the unit ball $(B_{x_0}(1), g)$. Thus, the function $f$ is proper and each level set is compact. Consequently, R. Haslhofer and R. M\"uller \cite{HM11comp} observe that there is a point where $f$ attains its infimum. By picking such point as the reference point, they derive the following growth estimate
\begin{equation}
	\label{estimatef2}
\frac{1}{4}(r(x)-5n)^2_{+}\leq f(x)+C_1 \leq \frac{1}{4}(r(x)+\sqrt{2n})^2, 
\end{equation}
	for $a_{+}:=\max\{0,a\}$, and $C_1:=\SS+|\nabla f|^2-f$ is a constant. 
	Furthermore, if $p_1$ and $p_2$ are two minimum points then 
	\[ d(p_1, p_2)\leq 5n+\sqrt{2n}.\]
Hence it follows that,
\begin{align*}
	0\leq |\nabla f|^2 &\leq \frac{1}{4}(r(x)+\sqrt{2n})^2,\\
	-\frac{n}{2}\leq -\Delta f &\leq -\frac{n}{2}+\frac{1}{4}(r(x)+\sqrt{2n})^2.
\end{align*}
For comparison of the distance function, one can rewrite $\rho=2\sqrt{f}$ then 
\begin{align*}
|\nabla\rho| &=|\frac{\nabla f}{\sqrt{f}}|\leq 1,\\
 \Delta \rho &\leq \frac{\Delta f}{\sqrt{f}}\leq \frac{n}{\rho}\approx \frac{c}{r}.
\end{align*}

For $\lambda\leq 0$, there is not much known how to bound the potential function without additional assumption. Here is one for a nontrivial expanding GRS with $\lambda=-\frac{1}{2}$, O. Munteanu  and J. Wang \cite{MW12weighted} prove that, for $r>2=-4\lambda$,
\begin{equation}
	\label{estimatef3}
	\frac{1}{4}r^2 - C r^{3/2}\sqrt{\ln{r}} \leq \sup_{\partial B_{x_0}(r)}(-f)(x)\leq \frac{1}{4}r^2+Cr
\end{equation} for some constant C.
%%%%%%%%%%%%%%%%%%%%%%%%%%%%%%
\subsection{\textbf{Estimates on Curvatures}}%%%%%%%%
%%%%%%%%%%%%%%%%%%%%%%%%%%%%%%
In this subsection we  examine the behavior of various curvatures. We first consider the case of  a shrinking GRS ($\lambda>0$), with the normalization  $$\int_M (4\pi)^{-n/2}e^{-f} dV =1.$$ In \cite{HM11comp}, R. Haslhofer and R. Muller first show that there is a quadratic upper bound for the scalar curvature 
\begin{align*}
	0\leq \SS(x)&\leq \frac{1}{4} (r(x)+\sqrt{2n})^2.
\end{align*}

Obtaining bounds over other components of the Riemannian curvature is more challenging. O. Munteanu and N. Sesum \cite{munse09} have derived the following integral curvature estimate, for any $\epsilon>0$, $$\int_M |\Rc|^2 e^{-\epsilon f}<\infty.$$ 
In particular, when $\epsilon=1$, it follows from (\ref{lapS}) that,
$$\int_M |\Rc|^2 e^{-f}=\frac{1}{2}\int_M \SS e^{-f}<\infty .$$
These estimates play an important roles in the classification of locally conformally flat ($\WW=0$) shrinking GRS or those with harmonic Weyl tensor ($\delta \WW=0$). 

%If certain assumption is imposed on a curvature term, it is possible to obtain stronger results. 

Recalling from the characterizing equation (\ref{grs}), one can expect that $\Rc$ (or Bakry-Emery Ricci) will play a crucial role in bounding its geometry of GRS. It is observed by O. Munteanu and M.-T. Wang \cite{MW11cur} that if Ricci curvature is bounded by $\ell$, then there exists a constant $a=a(n, \ell)>0$, such that the growth of the Riemannian curvature tensor is bounded by
\[|\RR|(x)\leq C (r(x)+1)^{a}.\] As a consequence of the above result, they prove that if the Ricci curvature is small enough then the soliton must be isometric to the Gaussian shrinker with $\lambda=\frac{1}{2}$.
Furthermore,  O. Munteanu and J. Wang \cite{MW17conical} prove that if the Ricci curvature is approaching zero at infinity, then for $p\geq 8n$, $a\leq \frac14 p$, 
\[\int_M |\Rc|^p f^a <\infty ,\] and this shrinking GRS must be smoothly asymptotic to a cone. B. Chow, P. Lu and B. Yang \cite{CLY2011} show that, outside a compact set, scalar curvature is bounded from below by $$\SS \geq \frac{c}{r(x)^2}$$ for some constant $c$.

It turns out that, in dimension four, bounded scalar curvature is sufficient to bound other components of the Riemannian curvature tensor. In \cite{MW15curv4d}, O. Munteanu and J. Wang show that if the scalar curvature $\SS$ is bounded, then there exist some constants $C_1$ and $C_2$, such that
\begin{align*}
	|\RR| &\leq C_1\SS, \\
	\RR &\geq -(\frac{C_2}{\ln(r+1)})^{1/4}.
\end{align*}
The second inequality implies that practically $\RR$ is non-negative at infinity. See also \cite{CRZ21} for an extension of this result.   Moreover, if the GSRS is a singularity model, then its curvature grows at most quadratically; while if the steady GRS is a singularity model, then it must have bounded curvature \cite{CFSZ2020}. 

 For the expanding GRS case ($\lambda<0$),  in dimension four, H.-D. Cao and T. Liu \cite{CL22expand} consider a complete expanding gradient  Ricci solitons with nonnegative Ricci curvature outside a compact set K. They prove that the Riemannian curvature tensor and its covariant derivative are bounded by its scalar curvature, i.e., for any $0\leq a<1$, on $M\backslash K$, there exists a constant $C_a$ such that
\begin{align*}
	|\RR|\leq C_a \SS^a, \\
	|\nabla \RR|\leq C_a \SS^a.
\end{align*}
Hence if the scalar curvature has at most polynomial decay at infinity, then  $|\RR|\leq C\SS$.  Furthermore, any 4-dimensional complete non-compact expanding gradient  Ricci soliton with non-negative Ricci curvature has finite asymptotic curvature ratio, provided that it has finite asymptotic scalar curvature ratio.

For expanders of dimension $n\geq 5$ and $\Rc\geq 0$,   H.-D. Cao, T. Liu, and J. Xie \cite{CLX23} prove that if the asymptotic scalar curvature ratio is finite, $$\limsup_{r \rightarrow \infty} \SS r^2< \infty,$$ then the Riemanian curvature tensor must have at least sub-quadratic decay, i.e., for any $0<\alpha<2$,  
\[  \limsup_{r\rightarrow \infty}|\RR| r^\alpha< \infty.\]
P.-Y. Chan \cite{Chan2023} shows that on a non-flat complete non-comppact expnder with $\SS \geq 0$, $\SS+|\nabla f|^2=-f$ and $\lim_{r\rightarrow \infty} f=-\infty$, there exist a positive constant $c$ such that $$\SS \geq c(\frac{n}{2}-f)^{1-\frac{n}{2}} e^{f-\frac{n}{2}}.$$ In the same paper, Chan also derives a lower bound for the scalar curvature under a different setting (non-positive Ricci).

In a systematic approach, H. Li, Y. Li, and B. Wang \cite{LLW21structure} develop a structure theory for non-collapsed Ricci shrinkers without any curvature condition. As an application, they are able to obtain curvature estimates that depend only on the non-collapsing constant and dimension.
% Cao-Zhou, 10 (Lemma 3.1): $\lambda>0$, the average scalar curvature over $D(r)$ is bounded above by $\frac{n}{2}$ ($n\lambda$). $D(r)=\{x, ~ 2\sqrt{f}(x)\leq r\}$ (here the normalization is $\Delta_f f+2\lambda f=n\lambda$). 
 
For the steady case ($\lambda=0$), with the convention $\SS+|\nabla f|^2=1$, assuming that $f\leq 0$ and $\lim_{r\rightarrow \infty} f=-\infty$, Chow, Lu and Yang \cite{CLY2011} provide a sharp lower bound for the scalar curvature $$\SS\geq \frac{2}{\sqrt{2n}+4}e^f.$$

\subsection{\textbf{Estimates on Volume}}
For a complete non-compact shrinking GRS  ($\lambda >0$) $M$, H.-D. Cao and D. Zhou \cite{caozhou10} show that it can have at most  Euclidean volume growth, i.e.,  for any $x_0\in M$ and sufficiently large $r$, there exists a constant $c$ such that,
	\[\text{Vol}(B_{x_0}(r))\leq c r^n.
	\]  
\begin{remark}
	The upper bound's power is optimal due to the construction of \cite{FIK03Kahler}. 
\end{remark}
\begin{remark}
	The above constant $c$ depends on the geometry (metric $g$) of the shrinking GRS. A crucial improvement was obtained in \cite{HM11comp} by showing that the constant $c$ can be chosen depending only on the dimension $n$. See also \cite{CRZ22} for related results. 
	
	 %(Normalize $\Delta_f f=2\lambda f+ c(1))$)
\end{remark}
%Corollary: finite integral of function with reasonable growth rate and $e^{-f}$-weighted volume.  \\
%Lower bound: must be infinite \cite{caozhu}. 

 O. Munteanu and J. Wang \cite{MW12weighted} prove a Calabi-Yau type lower bound,  for any $x_0 \in M$ and radius $r$, 
\[\text{Vol}(B_{x_0}(r))\geq c r. 
	\]

For  the steady case ($\lambda=0$), O. Munteanu and N. Sesum \cite{munse09} derive the following volume estimates: there exist positive constants $c, r_0, a$ such that, for any $r>r_0$,
	\[ce^{a\sqrt{r}}\geq \text{Vol}(B_{x_0}(r))\geq c^{-1} r.
	\]

In the expanding case ($\lambda<0$), an upper bound has been obtained by O. Munteanu and J. Wang \cite{MW14ADV}:
\[\text{Vol}(B_{x_0}(r))\leq C e^{\sqrt{n-1}r}.
	\]  
 If the scalar curvature $S\geq -\beta$ with $\beta\geq 0$, then J. Carrillo and L. Ni \cite{CN09CAG} prove that for some $r_0$ and any $r\geq r_0$, 
 \[\text{Vol}(B_{x_0}(r))\geq \text{Vol}(B_{x_0}(r_0)) (\frac{r+\alpha}{r_0+\alpha})^{n-2\beta},
	\]  where $\alpha=\sqrt{\beta-f(x_0)}$.

For a more detailed discussion on the volume estimate, please see P.-Y. Chan, Z. Ma and Y. Zhang \cite{CMZ2022} and references therein.

\subsection{Ends at Infinity}
It is an interesting question to count the number of ends for a complete non-compact GRS. For the expanding case ($\lambda<0$), if $\SS\geq  -\frac{n-1}{2}$, O. Munteanu and J. Wang \cite{MW12weighted} prove that  the GRS is either connected at infinity or isometric to $N \times \mathbb{R}$ for  a compact Einstein manifold $N$.

The K\"{a}hlerity assumption is particularly useful in this direction. For $\lambda>0 $, O. Munteanu and J. Wang \cite{MW15topo} prove that a K\"{a}hler GRS is connected at infinity, i.e., it has only one end. This is also true for any expanding Kähler GRS ($\lambda<0$) with proper potential function. For $\lambda=0$, O. Munteanu and N. Sesum \cite{munse09} study a K\"{a}hler GRS with Ricci curvature bounded below and a uniform non-collapsing condition. They conclude that the manifold is either connected at infinity (one end) or isometrically it is $M=\mathbb{R}\times N$ for a compact Ricci flat manifold $N$. Notice that the same conclusion is proved by O. Munteanu and J. Wang \cite{MW11mms, MW14ADV} only assuming completeness.

%Then, the authors obtain a sharp pointwise lower bound for the weight function of any smooth metric measure space, in terms of a lower bound on the associated Bakry-Émery curvature.

%: The authors continue here the study initiated by them in [Comm. Anal. Geom. 19 (2011), no. 3, 451–486; MR2843238; Comm. Anal. Geom. 20 (2012), no. 1, 55–94; MR2903101] of the geometry and topology of manifolds with Bakry-Émery curvature bounded from below, focusing on gradient Ricci solitons.

\subsection{Cone and Cylinder}
It has been observed that shrinking and steady gradient Ricci solitons have many common properties with manifolds of non-negative Ricci curvature. Consequently, a fundamental theme is to show that such a GRS is asymptotically close to a cone or a cylinder under certain curvature assumptions.  B. Kotschwar and L. Wang \cite{KW15rigid} show that two complete shrinking gradient Ricci solitons which are asymptotic to the same regular cone must be isometric to each other. This result essentially reduces the study of shrinking gradient Ricci solitons to their asymptotic cones. Similar ideas have been used in the study of self-shrinkers for the mean curvature flow \cite{wang14uniqueMCF}. 

%: In this paper, the authors prove an interesting rigidity theorem for asymptotically conical shrinking gradient Ricci solitons:
%Theorem 1.2. Suppose that two shrinking gradient Ricci solitons are asymptotic to the regular cone (E0,gc) along their ends respectively. Then there exist ends and a diffeomorphism.
%It has important corollaries:
%The key idea of the proof to Theorem 1.2 is to reduce the problem to a problem of backward uniqueness for the corresponding Ricci flow, which is a weakly parabolic system, and then to use "Carleman-type'' inequalities to prove the backward uniqueness, following [L. Escauriaza, G. A. Seregin and V. Šverák, Arch. Ration. Mech. Anal. 169 (2003), no. 2, 147–157; MR2005639].
%To derive these Carleman inequalities, the authors use a prolonged PDE-ODE system from previous work of the first author [Int. Math. Res. Not. IMRN 2010, no. 21, 4064–4097; MR2738351]. This is the technical heart of the paper under review.
%Note that similar ideas have also been used by the second author to study the self-shrinkers for the mean curvature flow [J. Amer. Math. Soc. 27 (2014), no. 3, 613–638; MR3194490].

Following the work of  Kotschwar-Wang, B. Chow and P. Lu \cite{CL15uniquecone} show that a complete non-compact shrinking GRS has at most one tangent cone at infinity, provided certain Ricci curvature decays.  The uniqueness is known to fail without a soliton structure (see references therein). %The authors use this to help establish that the canonical solution to Ricci flow determined by (M,g,f) converges smoothly outside a compact set (as t approaches the singular time) to an isometric copy of the complement of a compact set in a regular cone. This is then used directly to show that any two asymptotic cones are isometric.
Furthermore, they also provide a sufficient condition for the decay of the Ricci curvature: the decay of the scalar curvature, the boundedness of the full Riemannian curvature, and a non-collapsing condition.

Certain curvature conditions could be used to detect asymptotic cones. In \cite{MW17conical}, O. Munteanu and J. Wang prove that a shrinking GRS must be smoothly asymptotic to a cone if its Ricci curvature goes to zero at infinity. In dimension four, the assumption can be weakened significantly. In another article \cite{MW15curv4d}, they observed that it suffices to impose that the scalar curvature converges to zero at infinity.  %there exists a cone $E_0$ such that (M,g) is $C^k$ asymptotic to $E_0$ for all k. Also, they prove that for a compact gradient shrinking Ricci soliton of dimension n, the diameter of (M,g) has an upper bound of the form $diam(M)\leq c(n,inj(M))$, where inj(M) is the injectivity radius of (M,g).

For  the expanding case ($\lambda< 0$), C.-W. Chen and A. Deruelle \cite{CD15structureexpanding}  study the geometry at infinity of a complete expanding GRS. They show that the existence and uniqueness of an asymptotic cone is guaranteed by a finite asymptotic curvature ratio. In particular, it follows from non-negative curvature. In dimension four, H.-D. Cao and T. Liu \cite{CL22expand}, using new curvature estimates, prove the existence of a cone assuming nonnegative Ricci curvature and finite asymptotic scalar curvature ratio.

The cylindrical model also received tremendous interest. For the shrinking case ($\lambda>0$) with bounded curvature, O. Munteanu and J. Wang \cite{MW19structure} show that,  if the round cylinder $\mathbb{R}\times \mathbb{S}^{n-1}/\Gamma$ occurs as a limit for a sequence of points going to infinity along one end, then that end is asymptotic at infinity to the same round cylinder. In dimension four, if the scalar curvature is bounded below by a positive constant, then along each end  the soliton is asymptotic to a quotient of $\mathbb{R}\times \mathbb{S}^3$ or converges to a quotient of $\mathbb{R}^2\times \mathbb{S}^2$ along each integral curve of $\nabla f$ (the case that scalar curvature approaching zero at infinity must be smoothly asymptotic to a cone was proved in \cite{MW15curv4d}). 

%This result is applied to obtain structural results at infinity for four dimensional gradient shrinking Ricci solitons. It was previously known that such solitons with scalar curvature approaching zero at infinity must be smoothly asymptotic to a cone. For the case that the scalar curvature is bounded from below by a positive constant, we conclude that along each end

B. Kotschwar and L. Wang have studied the situation when a GRS is asymptotically similar to a cylinder. In \cite{KW17uniqueness}, they show that a shrinking GRS which agrees to infinite order at spatial infinity with one of the standard cylindrical metrics on $\mathbb{S}^k\times \RR^{n-k}$ ($k\geq 2$) along some end must be isometric to the cylinder on that end. If it's complete, then the GRS must be globally isometric  to either the cylinder or  its $Z_2$-quotient (when $k=n-1$) .

The splitting results are similar in spirit to those of cylindrical models. With the assumption of non-collapsing and bounded curvature, B. Chow and P. Lu  \cite{CL16split} conclude that such a sequence of complete-noncompact shrinking GRS must sub-converge (in the Cheeger-Gromov sense) to a limit which splits a line. 
For a steady 4-dimensional GRS, B. Chow, Y. Deng and Z. Ma \cite{CDM22dimensionreduction} prove that if the soliton dimension reduces to 3-manifolds then it is either strongly dimension reduces to a quotient of $\mathbb{S}^3$, or weakly dimension reduces to the Bryant 3-soliton. Here by dimension reduction we  mean a rescaled sequence of regions on the manifold converges  to a product with a line factor in the pointed Cheeger-Gromov sense. 

%This paper is about the asymptotic geometry of 4-dimensional steady gradient Ricci solitons and 4-dimensional steady gradient Ricci soliton singularity models. 
% A manifold $(M^n,g,f)$ dimension reduces to $(n-1)$-manifolds if a subsequence of metrics converges to $N^{n-1}\times R$ in the $C^\infty$ pointed Cheeger-Gromov sense, where $(N,g_N(t))$ is an $(n-1)$-dimensional complete ancient Ricci flow with bounded curvature.

% When the steady Ricci solitons are Kähler-Ricci solitons, they classify steady gradient Kähler-Ricci soliton singularity models of complex dimension 2 (Theorem 1.7).
%The authors prove that a 4-dimensional steady gradient Ricci soliton whose dimension reduces to  (Theorem 1.3). Its proof is related to linear curvature decay. 
%\begin{itemize}
%	\item \cite{munse09}. Note the fact, which motivates people to study the issue of existence of harmonic functions on these manifolds. Moreover, it is known that the existence of certain classes of harmonic functions is related to the existence of ends of the manifold. Thus in the second part of this paper, the authors address the issue of existence of harmonic functions on gradient shrinking Kähler and gradient steady Ricci solitons. They firstly prove that for a gradient steady Ricci soliton, if u is harmonic with $\int_M |\nabla u|^2< \infty$, then u is constant on M. As a consequence, the authors derive that any gradient steady Ricci soliton has at most one nonparabolic end. Based on this, consequences to the structure of shrinking and steady solitons at infinity are also discussed.

If a shrinking GRS is asymptotic to a Kähler cone along one end $V$, B. Kotschwar \cite{Kotschwar18Kahlercone} proves that such a gradient Ricci soliton is Kähler on some neighborhood of infinity of that end $V$, and it is globally Kähler if the soliton is complete.

%	\item \cite{BCDMZ22}: In this paper we consider 4-dimensional steady soliton singularity models, i.e., complete steady gradient Ricci solitons that arise as the rescaled limit of a finite time singular solution of the Ricci flow on a closed 4-manifold. In particular, we study the geometry at infinity of such Ricci solitons under the assumption that their tangent flow at infinity is the product of R with a 3-dimensional spherical space form. We also classify the tangent flows at infinity of 4-dimensional steady soliton singularity models in general.

%\end{itemize}

%%%%%%%%%%%%%%%%%%%%%%%%%%%%%%%%%%%%
%\input{class0202}
%\input{class0302}
\section{\textbf{Classification and Rigidity}}
%06112024

In this section, we review classification results on the subject. Due to smoothing effects of the Ricci flow equation, a GRS being an equilibrium is expected to be rich in symmetry, particularly in low dimensions. The literature is consistent with this principle. First, in dimensions two and three, their geometry is relatively better understood (except the expanding case in dimension three). In higher dimensions, the situation is significantly more complicated mainly due to the non-triviality of the Weyl tensor ($\WW$), which is vacuously zero for dimension less than four. Thus, it is reasonable to impose certain  restriction on the curvature operator or one of its components: the Weyl, the Ricci, or the scalar curvature. Dimension four also exhibit distinguishing phenomenon and there have been important advancements accordingly.  
%imposing certain restrictions on the curvature operator. For example, a natural development is to impose conditions on that Weyl tensor in higher dimensions. 

\subsection{Dimensions Two and Three}
For $n=2$, any GRS on a closed surface must have constant curvature. For open and non-compact surfaces, there is a non-trivial Killing vector field indicating rotational symmetry \cite{Hsurface, bm15}. For a complete list, we refer to \cite{chow2023}, and references therein.

For $n=3$, there is the Hamilton-Ivey estimate for the Ricci flow equation. It effectively says that approaching a finite singular time, the most negative sectional curvature is relatively small compared to the most positive one. Consequently, via a careful point-picking procedure, one can obtain a non-negatively curved limit. It has amplified effects on the classification scheme. G. Perelman \cite{perelman2} first shows that a complete non-flat, non-collapsed shrinking GRS with non-negative sectional curvature must be either $\mathbb{S}^3$ with the round metric, the round infinite cylinder $\mathbb{S}^2\times \mathbb{R}$, or their quotients. Later, L. Ni and N. Wallach \cite{niwa08}, A. Naber \cite{naber07}, and H.-D. Cao, B.-L. Chen, X.-P. Zhu \cite{CCZ08} succeed in removing the non-collasing and all curvature assumptions. 

%being $\kappa$-non-collapsed on all scales for a $\kappa>0$

Similarly, S. Brendle proves that a complete non-flat non-collapsed steady GRS must be rotationally symmetric and, consequently, isometric to the Bryant soliton \cite{b12rot}. The non-collapsing assumption is crucial in his proof. In general, all $3$-dimensional steady GRS are non-negatively curved and each is asymptotic to a sector of angle $\alpha\in [0, \pi]$ \cite{chenbl09}. The Bryant soliton corresponds to $\alpha=0$ while the product of a $2$-dimensional cigar with a line corresponds to $\alpha=\pi$. Recently, Y. Lai shows that such a solition must have $O(2)$ symmetry and, furthermore, she constructs examples (called flying wings) such that $\alpha$ assumes an arbitrary value in $[0, \pi]$ \cite{Lai_fylingwings, Lai_O2}.    

On the other hand, the expanding case is significantly less rigid, there is no general classification to our best knowledge. P. Petersen and W. Wylie \cite{pewy10} prove that the only $3$-dimensional expanding GRS with constant scalar curvature are  quotients of $\mathbb{R}^3$, $\mathbb{H}^2\times \mathbb{R}$, or $\mathbb{H}^3$. 

%%%%%%%%%%%%%%%%%%%%
\subsection{Higher Dimensions} %%%%
%%%%%%%%%%%%%%%%
It turns out that having constant scalar curvature is a restrictive condition (this is automatically true for an Einstein metric). In the steady case, constant scalar curvature immediately implies that $\Rc\equiv 0$ \cite{PW09rigid}. Indeed, P. Petersen and W. Wylie \cite{PW09rigid} have obtained a stronger result only assuming that scalar curvature $\SS$ attains its minimum value. %Furthermore, M. Eminenti, G. La Nave and C. Mantegazza \cite{ELM08} prove that a compact shrinking GRS with constant Ricci curvature and zero Weyl tensor is rigid if and only if its scalar curvature is constant. 

%utilizing the Hamilton-Ivey estimate \cite{ivey}, shrinking GRS must either have positive sectional curvature hence compact, or it must be a quotient of a product of $\RR$ and a 2-dimensional shrinker. Perelman \cite{perelman2} proved an analogous theorem. {\textcolor{red} {(A little more details? HT I'm not sure which one)}%	Complete non-compact steady GRS must be Bryant soliton. {\textcolor{red} {(Gaussian shrinking GRS on $\RR^3$, the round cylinder shrinker on $S\sp 2 \times \RR$)}%{\textcolor{red} {A steady GRS is flat or rotationally symmetric (that is, a Bryant Soliton) by \cite{caochen12}. where shall we put this?? \cite{PW09rigid} All homogeneous gradient Ricci solitons are rigid.}		%Hamilton \cite{Hsurface}, Bryant \cite{}, Ivey \cite{}, Berstein and Metter \cite{bm15} (for a complete list, we refer to \cite{chow2023}, and references therein), completely classified GRS on closed surfaces must be either the round sphere, projective space, flat torus, Klein bottle, or a surface of constant negative curvature with genus at least 2. Complete non-compact steady or shrinking GRS can only be the Cigar soliton or flat Euclidean space. {\textcolor{red} {(how to describe the expanding case?, do we need to? HT: In dimension two, we can refer to \cite{bm15}) }

Using the theory of isoparametric functions, M. Fernandez-Lopez and E. Garcia-Rio \cite{FG16csc} prove that if a $n$-dimensional GRS has constant scalar then $\SS=k\lambda$ for some integer $0\leq k \leq n$. Furthermore, they show that such a soliton is rigid if it's one of the following: 1)  $k=n-1$; 2) its Ricci tensor has constant rank; 3) or it's  a K\"{a}hler GRS in real dimension four and six.  %n case $\lambda>0$, it is observed that $k\neq 1$. 

The Weyl curvature tensor also plays an important role in understanding the geometry of GRS. The Weyl tensor measures locally conformal flatness, and a Riemannian manifold has vanishing Weyl tensor if it's rotationally symmetric. By analyzing the corresponding ODEs, B. Kotschwar \cite{kb08} gives a classification of all rotationally symmetric gradient Ricci solitons with given diffeomorphism types of  $\mathbb{S}^n$, $\RR^n$, or $S^{n-1}\times \RR$. It follows from the works of \cite{niwa08, zhangzh09, cw07, pewy10, cm11} that if $\WW\equiv 0$, then there is no other complete shrinking GRS besides finite quotients of those above. In the steady case, H.-D. Cao and Q. Chen show that a complete non-compact locally conformally flat gradient Ricci soliton must be either flat or isometric to the Bryant soliton \cite{caochen12}. Z.-H. Zhang \cite{zhangzh09} made an important observation  that locally conformal flatness implies non-negative curvature operator.  For the expanding case, to our best knowledge, there is no analogous results in arbitrary dimension  yet other than in dimension four.

The condition of vanishing Weyl tensor $\WW\equiv 0$ can be weakened to harmonic Weyl tensor $\delta\WW\equiv 0$ \cite{fega11, munse09} without introducing any new soliton model. In a similar approach, there are other classification results based on vanishing of the Bach tensor \cite{caochen11}, a fourth-order vanishing condition on the Weyl tensor \cite{CMM17vanishing}, or assumptions on the radial sectional curvature \cite{pewy10}. 

%In this case, a complete shrinking GRS with harmonic Weyl tensor must be a finite quotient of $\mathbb{R}^n$, $\mathbb{S}^{n-1}\times \mathbb{R}$, or  $\mathbb{S}^{n}$. \\

%\cite{CMM17vanishing}: We classify complete gradient Ricci solitons satisfying  improving previously known results. More precisely, we show that any n-dimensional ($n\geq 4$) gradient shrinking Ricci soliton with fourth order divergence-free Weyl tensor is either Einstein, or a finite quotient of $N^{n-k}\times R^k$, ($k>0$), the product of a Einstein manifold $N^{n-k}$ with the Gaussian shrinking soliton $R^k$. The technique applies also to the steady and expanding cases in all dimensions. In particular, we prove that a three dimensional gradient steady soliton with third order divergence-free Cotton tensor, i.e. with vanishing double divergence of the Bach tensor, is either flat or isometric to the Bryant soliton.''
			
			% \cite{caotran1}: 

GRS with positive isotropic curvature (PIC) attracts significant interest in the study; this is a notion first introduced by M. Micallef and J. Moore \cite{MM88}. S. Brendle and R. Schoen also consider natural generalizations, PIC1 and PIC2 \cite{bs072} (generally speaking, PIC$k$ means the product metric with an Euclidean factor $\mathbb{R}^k$ is PIC). Since these conditions are preserved along the Ricci flow \cite{HPIC, bs072, nguyen10}, they persist to singularity models and gradient Ricci solitons.  For 4-dimensional Einstein manifolds, it is well-known that positive isotropic curvature is equivalent to 2-positive curvature operator. 

In higher dimensions, S. Brendle \cite{brendle10einstein} proves that closed Einstein manifolds with nonnegative isotropic curvature must be locally symmetric. When dimension $n\geq 12$, K. Naff \cite{naff2019shrinking} proves that a complete non-flat  shrinking gradient Ricci solitons with uniformly PIC must be a quotient of either the round sphere $\mathbb{S}^n$ or the cylinder $\mathbb{S}^{n-1}\times \mathbb{R}$. The proof relies crucially on the general theory developed by S. Brendle \cite{brendle19} on Ricci flows. Furthermore, a complete shrinking GRS that is strictly PIC and weakly PIC2 must be a quotient of either the round sphere $\mathbb{S}^n$ or the cylinder $\mathbb{S}^{n-1}\times \mathbb{R}$, provided that dimension $n\geq 5$. In dimension four, there is an improvement by \cite{caoxie2023}. Furthermore, in \cite{LN20}, X. Li and L. Ni give a classification with weakly PIC1 and it is easy to see among the examples only the sphere and cylinder and their quotients have strictly PIC.

% The result for n\geq 12 with uniform PIC relies on Brendle's work (MR3997128). Naff's result with strictly PIC and weakly PIC2 actually follows from Li-Ni (MR4098770), as Li-Ni gave a classification with wealy PIC1 and it is easy to see among the examples only the sphere and cylinder and their quotients have strictly PIC.
			
In a slightly different approach \cite{CMM15}, G. Catino, P. Mastrolia and D. Monticelli prove that for a complete 
steady GRS $(M^n, g_{ij} , f)$  with nonnegative sectional curvature and $$\liminf_{r\rightarrow \infty} \frac{1}{r} \int_{B_{r}(O)} R=0$$ (zero re-scaled average scalar curvature), then $M^n$ is isometric to a quotient of Gaussian soliton or $R^{n-2}\times \Sigma^2$, where $\Sigma^2$ is the cigar soliton. 
The only complete steady GRS with nonnegative sectional curvature and less than quadratic volume growth are quotients of $R^{n-2}\times \Sigma^2$.

In the expanding case, if $(M^n, g_{ij} , f)$ is a complete expanding GRS
with nonnegative sectional curvature and $S \in L^1(M)$, then M is isometric to
a quotient of the Gaussian soliton $R^n$ \cite{CMM15}. If $|\nabla f| \in  L^p(M, e^{-f}dvol)$ for some $1 \leq p \leq \infty$, then the expanding GRS $(M^n, g_{ij} , f)$ is trivial \cite{prs09}. Consequently, if scalar curvature $S \geq 0$ and $S \in L^1(M, e^{-f}dvol)$, then M is isometric
to the standard Euclidean space.

P. Petersen and W. Wylie \cite{PW09rigid} prove that a shrinking (expanding) gradient soliton is rigid (in the sense that is a flat bundle of $N\times\mathbb{R}^k$ with $N$ being Einstein) if and only if it has
constant scalar curvature and is radially flat. Furthermore, all complete non-compact shrinking gradient solitons of co-homogeneity one with nonnegative Ricci curvature and $\text{K}(\cdot,\nabla f) \geq 0$ are rigid	\cite{PW09symmetry}.	
	
\subsection{Dimension Four} 
Dimension four is unique due to, among others, a certain decomposition of the Lie algebra governing a natural principal bundle over a Riemannian manifold. Consequently, there have been several important improvements and developments. 

A crucial analogue of Hamilton-Perelman results is obtained by A. Naber proving that a $4$-dimensional complete non-compact shrinking GRS with bounded non-negative curvature operator must be a finite quotient of $\RR^4$,  $\mathbb{S}^2\times \RR^2$ or $\mathbb{S}^3\times \RR$ \cite{naber07}. The result paves the way for significant later advancement among which is the classification of K\"{a}hler shrinker surfaces \cite{tian1, WZ04toric, CS2018classification, BCCD22KahlerRicci, LW23}. 

Considering a shrinking GRS with constant scalar curvature, X. Cheng and D. Zhou \cite{CZ2023rigidity} extend the work initiated by M. Fernandez-Lopez and E. Garcia-Rio \cite{FG16csc}.  They prove that a four-dimensional complete shrinking GRS is rigid  if and only if it has constant scalar curvature. %, that is, it is either Einstein, or a finite quotient of a Gaussian shrinking soliton $\RR^4$,  $\mathbb{S}^2\times \RR^2$ or $\mathbb{S}^3\times \RR$.  

In dimension four, the Weyl curvature tensor  is naturally decomposed into self-dual and anti-self-dual parts $\WW^{\pm}$ (depending on a choice of orientation). A Bochner-Weitzenb\"ock type formula is derived by the authors in \cite{caotran1}. Consequently, it leads to several classification and rigidity results based on conditions on the Weyl tensor and the potential functions. The vanishing of the self-dual part $\WW^+$ has been studied by \cite{cw11}. It is interesting to note that such a consideration would rule out the case of $\mathbb{S}^2\times \RR^2$ while still adimitting $\mathbb{CP}^2$. Consequently, it is later observed that it suffices to only impose the vanishing of higher order terms. In \cite{www2014, www20142}, J. Wu, P. Wu, and W. Wylie prove that a four-dimensional complete, oriented, non-Einstein shrinking GRS with harmonic self-dual Weyl tensor (the divergence of $\WW^{+}$ is vansihing), is necessarily isometric to a finite quotient of $\mathbb{S}^2\times \RR^2$ or $\mathbb{S}^3\times \RR$.

Analogous results have been obtained for steady and expanding solitons. J. Kim  \cite{kim17harmonic} shows that a 4-dimensional complete steady GRS with harmonic Weyl curvature is either Ricci flat or isometric to the Bryant soliton.  Furthermore, any 4-dimensional complete expanding gradient Ricci soliton with harmonic Weyl curvature must be one of the following: Einstein, rigid, or conformally flat. In \cite{CDM22dimensionreduction}, B. Chow, Y. Deng and Z. Ma classify steady four-dimensional GRS with non-negative Ricci curvature outside a compact set. They prove that 4-dimensional steady gradient Ricci soliton singularity models with nonnegative Ricci curvature outside a compact set must be either Ricci-flat ALE 4-manifolds or dimension reduce to 3-dimensional manifolds.

%They show that a shrinking GRS with bounded curvature must be a finite quotient of $\RR^4$, $S^{3}\times \RR$, $S^n$, or $CP^2$, and a steady GRS must be a Bryant soliton or flat.

%	\item : "The paper contains two main results.  Second, a four-dimensional complete non-Einstein gradient Kähler-Ricci soliton with constant scalar curvature must be isometric to a finite quotient of the Riemannian product of a surface and C. (The completeness hypothesis, obviously necessary in both cases, is left unstated.) The proofs involve constructing a weighted subharmonic function."
	
%	\item 
%\end{enumerate} 

 %$n=4$, GRS, constant scalar and $S\neq 2\lambda$ then rigid. If $S=2\lambda$ then $\Rc\geq 0$ (shrinking) and $\Rc\leq 0$ (expanding).
 
% . Petersen-Wylie \cite{PW09rigid}. 
%[Int. Math. Res. Not. IMRN 2008, no. 4, Art. ID rnm152; MR2424175] by removing the nonnegative curvature operator, a pinching condition and a curvature growth condition. 
%the topic of {K}\"{a}hler GRS already received tremendous interest: \cite{caohd96,  CZ12Kahler, FIK03Kahler, MW15topo, CD2020expanding, DZ2020rigidity, CS2018classification, CF16conical, BCCD22KahlerRicci}. 
%\item Complete classification (under review): Yu Li-Bing Wang \cite{LW23}

In \cite{LNW16fourPIC}, X. Li, L. Ni and K. Wang prove that any four-dimensional complete shrinking GRS with positive isotropic curvature is either a quotient of a round sphere $\mathbb{S}^4$ or the cylinder $\mathbb{S}^{3}\times \mathbb{R}$. This generalizes an earlier result of L. Ni and N. R. Wallach \cite{niwa08}. In addition, that paper also has a classification of 4-dimensional  complete shrinkers with nonnegative isotropic curvature, which generalizes the result of A. Naber \cite{naber07}. 

\subsection{K\"{a}hler Ricci Solitons}

As mentioned in the intruduction, the K\"{a}hler case  has extensive literature; see, for examples, \cite{tian1, WZ04toric, caohd09, MW15topo, CS2018classification, CF16conical, CD2020expanding, DZ2020rigidity} and references therein. In particular, tremendous recent developments led to the classification of all shrinking K\"{a}hler  Ricci solitons in complex dimension 2 \cite{CDSexpandshriking19, CCD22finite, BCCD22KahlerRicci, LW23}. Furthermore, for any sign of $\lambda$, it is shown that a KGRS must be toric under a generic assumption \cite{tran23contact, tran24toric}.

\subsection{Rigidity in Strong Sense}
In this subsection, we discuss a different type of ``Rigidity" from earlier sections, that is, if a GRS is close to a model in certain sense, then it must be exactly that model. This direction could be considered as a combination of the classification results and asymptotic analysis. In \cite{LW21rigidcylinder}, Y. Li and B. Wang show that the round cylinders are rigid in the space of Ricci shrinkers, i.e., any Ricci shrinker that is sufficiently close to $\mathbb{S}^{n-1}\times \mathbb{R}$ in the pointed Gromov-Hausdorff topology must be isometric to $\mathbb{S}^{n-1}\times \mathbb{R}$. Moreover, T. Colding and W. Minicozzi \cite{colding2022singularities} show that cylindrical shrinkers $\mathbb{S}^{n-l}\times \mathbb{R}^{l}$ are {\it strongly rigid} for any $l$.

In the table below, we list a few notable results in this direction (disclaimer: it is by no means exhausted).

\begin{center}
	\begin{tabular}{ |c|c|c| } 
		\hline
		Assumptions & Models & Reference \\
		\hline 
		$\lambda>0$, $0\leq K< C$, decay of $\nabla \Rc$ & Rigid & \cite{cai15} \\ 
		\hline 
		$\lambda>0$, $\Rc\geq 0$, pinched $|W|$ & $\WW\equiv 0$ & \cite{catino13pinched} \\ 
		\hline
		$\lambda=0$, $K>0$, asymptotically cylindrical & Rotational Symmetry & \cite{brendle14rotahigh} \\ 
		\hline
		$\lambda<0$, $K>0$, asymptotically conical & Rotational Symmetry & \cite{chodosh14}\\ 
		\hline
	\end{tabular}
\end{center}

See also \cite{CRT22preprint} for some improvement on the result by \cite{catino13pinched}.

%\section{Int.}
%%%%%%%%%%%%%%%%%%%%%%%%%%%%%%%%%%%%%

%\def\cprime{$'$}
\bibliographystyle{plain}
\bibliography{biomorse}

%%%%%%%%%%%%%%%%%%%%%%%%%%%%%%%%
\end{document}